\documentstyle[11pt, amssymb]{article}

\begin{document}

\begin{flushleft}
{\Large NOTES ON MIXED TEICHM\"{U}LLER MOTIVES}
\end{flushleft}

\noindent
\hrulefill

\begin{flushleft}
{\large TAKASHI ICHIKAWA} 
\end{flushleft}
\vspace{2ex}

\noindent
{\bf Abstract} 

\noindent
As a higher genus version of universal mixed elliptic motives by Hain and Matsumoto, 
we consider mixed Teichm\"{u}ller motives as certain motivic local systems 
on the moduli space of pointed curves. 
We show that the category of mixed Teichm\"{u}ller motives is equivalent to 
a full subcategory of a certain product category of 
mixed Tate motives over ${\mathbb Z}$ and universal mixed elliptic motives. 
Furthermore, we show that unipotent fundamental torsors for the universal open curve 
give rise to a pro-object in the category of mixed Teichm\"{u}ller motives. 
Our results can give a realization of motivic correlators proposed by Goncharov. 
\vspace{4ex}

\noindent
{\bf 1. Introduction} 

\noindent
The aim of this paper is study mixed Teichm\"{u}ller motives 
defined as a higher genus version of mixed Tate motives (unramified) over ${\mathbb Z}$ 
(see Deligne-Goncharov \cite{DG} for the definitive reference) and 
universal mixed elliptic motives by Hain and Matsumoto \cite{HM}. 
For nonnegative integers $g, n$ such that $n, 2g - 2 + n > 0$, 
let ${\cal M}_{g, \vec{n}}$ denote the moduli stack over ${\mathbb Z}$ of 
proper smooth curves of genus $g$ with $n$ tangent vectors, 
and ${\cal M}_{g, \vec{n}}^{\rm an}$ denote the associated complex orbifold. 
Then we consider {\it mixed Teichm\"{u}ller motives} for ${\cal M}_{g, \vec{n}}$ 
as motivic unipotent local systems on ${\cal M}_{g, \vec{n}}^{\rm an}$ such that 
their pure parts are generated by the cohomology groups (with Tate twists) of corresponding curves, 
and that they give mixed Tate motives over ${\mathbb Z}$ on the points at infinity. 
From the viewpoint of (arithmetic) algebraic geometry, 
it seems interesting to study the relationship with sheaves of mixed motives 
over ${\cal M}_{g, \vec{n}}$ given by Voevodsky, Ayoub and others \cite{VSF, A}. 
Our standpoint based on number theory is to develop the theory of 
mixed Teichm\"{u}ller motives 
in order to study the arithmetic properties of periods and $L$-values of curves. 

One of main results of this paper describes the structure of 
the category of mixed Teichm\"{u}ller motives 
using the theory of Teichm\"{u}ller groupoid for ${\cal M}_{g, \vec{n}}$ 
defined as its fundamental groupoid with tangential base points at infinity 
(cf. \cite{BK1, BK2, FG, Gr2, HLS, I4, MS, NS}). 
We describe the category of mixed Teichm\"{u}ller motives in terms of 
Teichm\"{u}ller's Lego game proposed by Grothendieck \cite{Gr2}. 
Namely, we show that the category of mixed Teichm\"{u}ller motives is equivalent to 
a full subcategory of a certain product category of 
mixed Tate motives over ${\mathbb Z}$ and universal mixed elliptic motives. 
As its application, by works of Brown \cite{Br}, Terasoma \cite{T}, 
Deligne-Goncharov \cite{DG} and Hain-Matsumoto \cite{HM} 
on mixed Tate motives over ${\mathbb Z}$ and universal mixed elliptic motives, 
one can describe the motivic Galois group of the category of mixed Teichm\"{u}ller motives. 
 
Another main result gives important examples of mixed Teichm\"{u}ller motives. 
Denote by ${\cal C}$ the universal curve over ${\cal M}_{g, \vec{n}}$ 
with sections (and tangential vectors) $\sigma_{1}, ..., \sigma_{n}$, 
and put 
$$
{\cal C}^{\circ} =  {\cal C} - \bigcup_{k = 1}^{n} {\rm Im}(\sigma_{k}). 
$$
Then for $i, j \in \{ 1,..., n \}$, 
the sets $\pi_{1} \left( {\cal C}^{\circ}_{s}; (\sigma_{i})_{s},  (\sigma_{j})_{s} \right)$ 
$(s \in {\cal M}_{g, \vec{n}}^{\rm an})$ of homotopy classes of paths from  $(\sigma_{i})_{s}$ to $(\sigma_{j})_{s}$ give a motivic local system on 
${\cal M}_{g, \vec{n}}$ as a torsor over 
${\displaystyle \lim_{\longleftarrow} 
{\mathbb Q} \left[ \pi_{1} ({\cal C}^{\circ}_{s}) \right] / I^{m}}$, 
where $I$ denotes the augmentation ideal of 
${\mathbb Q} \left[ \pi_{1} ({\cal C}^{\circ}_{s}) \right]$.  
We show that this local system gives rise to a pro-object 
in the category of mixed Teichm\"{u}ller motives for ${\cal M}_{g, \vec{n}}$. 
Therefore, one has the associated motivic monodromy representation which 
substantially gives the motivic correlators proposed by Goncharov \cite{G} 
under the existence of the abelian category of mixed motives over fields. 
\vspace{4ex}

\noindent
{\bf 2. Generalized Tate curve}  
\vspace{2ex}

\noindent
{\it 2.1. Degenerate curve} 

\noindent
We review the well known correspondence 
between certain graphs and degenerate pointed curves, 
where a (pointed) curve is called {\it degenerate} if it is a stable (pointed) curve and 
the normalization of its irreducible components are all projective (pointed) lines. 
A {\it graph} $\Delta = (V, E, T)$ means a collection  
of 3 finite sets $V$ of vertices, $E$ of edges, $T$ of tails 
and 2 boundary maps 
$$
b : T \rightarrow V, 
\ \ b : E \longrightarrow \left( V \cup \{ \mbox{unordered pairs of elements of $V$} \} \right) 
$$
such that the geometric realization of $\Delta$ is connected. 
A graph $\Delta$ is called {\it stable} 
if its each vertex has degree $\geq 3$,  
i.e. has at least $3$ branches. 
Then for a degenerate pointed curve, 
its dual graph $\Delta = (V, E, T)$ is given by the correspondence: 
$$
\begin{array}{lcl}
V & \longleftrightarrow & 
\{ \mbox{irreducible components of the curve} \}, 
\\
E & \longleftrightarrow & 
\{ \mbox{singular points on the curve} \}, 
\\
T & \longleftrightarrow & 
\{ \mbox{marked points on the curve} \} 
\end{array}
$$
such that an edge (resp. a tail) of $\Delta$ has a vertex as its boundary 
if the corresponding singular (resp. marked) point belongs 
to the corresponding component. 
Denote by $\sharp X$ the number of elements of a finite set $X$. 
Under fixing a bijection 
$\nu : T \stackrel{\sim}{\rightarrow} \{ 1, ... , \sharp T \}$, 
which we call a numbering of $T$, 
a stable graph $\Delta = (V, E, T)$ becomes the dual graph 
of a degenerate $\sharp T$-pointed curve of genus 
${\rm rank}_{\mathbb Z} H_{1}(\Delta, {\mathbb Z})$ 
and that each tail $h \in T$ corresponds to the $\nu(h)$th marked point. 
In particular, a stable graph without tail is the dual graph of 
a degenerate (unpointed) curve by this correspondence. 
If $\Delta$ is trivalent, i.e. any vertex of $\Delta$ has just $3$ branches, 
then a degenerate $\sharp T$-pointed curve with dual graph $\Delta$ 
is maximally degenerate. 

An {\it orientation} of $\Delta = (V, E, T)$ means 
giving an orientation of each $e \in E$. 
Under an orientation of $\Delta$, 
denote by 
$$
\pm E 
\ = \ 
\{ e, -e \ | \ e \in E \} 
$$
the set of oriented edges, 
and by $v_{h}$ the terminal vertex of $h \in \pm E$ 
(resp. the boundary vertex $b(h)$ of $h \in T)$. 
If $e$ is a loop, then in fact $v_{e} = v_{-e}$. 
For each $h \in \pm E$, 
let $| h | \in E$ be the edge $h$ without orientation. 
\vspace{2ex}

\noindent
{\it 2.2. Universal Schottky group} 

\noindent
Let $\Delta = (V, E)$ be a stable graph without tail. 
Fix an orientation of $\Delta$, 
and take a subset ${\cal E}$ of $\pm E$ 
whose complement ${\cal E}_{\infty}$ satisfies the condition that 
$$
{\cal E}_{\infty} 
\cap 
\{ -h \ | \ h \in {\cal E}_{\infty} \} 
\ = \ 
\emptyset, 
$$ 
and that $v_{h} \neq v_{h'}$ for any distinct $h, h' \in {\cal E}_{\infty}$. 
We attach variables $\alpha_{h}$ for $h \in {\cal E}$ 
and $q_{e} = q_{-e}$ for $e \in E$. 
Let $A_{0}$ be the ${\mathbb Z}$-algebra generated by 
$\alpha_{h}$ $(h \in {\cal E})$, 
$1/(\alpha_{e} - \alpha_{-e})$ $(e, -e \in {\cal E})$ 
and $1/(\alpha_{h} - \alpha_{h'})$ 
$(h, h' \in {\cal E}$ with $h \neq h'$ and $v_{h} = v_{h'})$, 
and let 
$$ 
A \ = \ A_{0} [[q_{e} \ (e \in E)]], \ \ 
B \ = \ A \left[ \prod_{e \in E} q_{e}^{-1} \right]. 
$$ 
According to \cite[Section 2]{I2}, 
we construct the universal Schottky group $\Gamma$ 
associated with oriented $\Delta$ and ${\cal E}$ as follows. 
For $h \in \pm E$, 
let $\phi_{h}$ be the element of $PGL_{2}(B) = GL_{2}(B)/B^{\times}$ given by 
$$
\phi_{h} \ = \ 
\frac{1}{\alpha_{h} - \alpha_{-h}} 
\left( \begin{array}{cc} 
\alpha_{h} - \alpha_{-h} q_{h} & - \alpha_{h} \alpha_{-h} (1 - q_{h}) 
\\ 1 - q_{h} & -\alpha_{-h} + \alpha_{h} q_{h} 
\end{array} \right) 
\ {\rm mod}(B^{\times}), 
$$
where $\alpha_{h}$ (resp. $\alpha_{-h})$ means $\infty$ 
if $h$ (resp. $-h)$ belongs to ${\cal E}_{\infty}$. 
Then 
$$
\frac{\phi_{h}(z) - \alpha_{h}}{z - \alpha_{h}} 
\ = \ 
q_{h} \frac{\phi_{h}(z) - \alpha_{-h}}{z - \alpha_{-h}} 
\ \ (z \in {\mathbb P}^{1}), 
$$
where $PGL_{2}$ acts on ${\mathbb P}^{1}$ by linear fractional transformation. 
\vspace{2ex}

\noindent
{\it 2.3. Generalized Tate curve} 

\noindent
For any reduced path $\rho = h(1) \cdot h(2) \cdots h(l)$ 
which is the product of oriented edges $h(1), ... ,h(l)$ 
such that $h(i) \neq -h(i+1)$ and $v_{h(i)} = v_{-h(i+1)}$, 
one can associate an element $\rho^{*}$ of $PGL_{2}(B)$ 
having reduced expression 
$\phi_{h(l)} \phi_{h(l-1)} \cdots \phi_{h(1)}$. 
Fix a base point $v_{b}$ of $V$, 
and consider the fundamental group 
$\pi_{1} (\Delta, v_{b})$ which is a free group 
of rank $g = {\rm rank}_{\mathbb Z} H_{1}(\Delta, {\mathbb Z})$. 
Then the correspondence $\rho \mapsto \rho^{*}$ 
gives an injective anti-homomorphism 
$\pi_{1} (\Delta, v_{b}) \rightarrow PGL_{2}(B)$ 
whose image is denoted by $\Gamma$. 
As is shown in \cite[Proposition 1.3]{I2}, 
each $\gamma \in \Gamma$ has 
its attractive (resp. repulsive) fixed point $a$ (resp. $a')$ 
in ${\mathbb P}^{1}(\Omega)$ 
and its multiplier $b \in \sum_{e \in E} A q_{e}$ 
which satisfy that 
$$
\frac{\gamma(z) - a}{z - a} 
\ = \ 
b \frac{\gamma(z) - a'}{z - a'} 
\ \ \ 
(z \in {\mathbb P}^{1}(\Omega)). 
$$ 

It is shown in \cite[Section 3]{I2} 
(and had been shown in \cite[Section 2]{IN} when $\Delta$ is trivalent and has no loop) 
that for any stable graph $\Delta = (V, E)$ without tail, 
there exists a stable curve $C_{\Delta}$ of genus $g$ over $A$ 
which satisfies the following: 

\begin{itemize}

\item 
The closed fiber $C_{\Delta} \otimes_{A} A_{0}$ of $C_{\Delta}$ 
obtained by substituting $q_{e} = 0$ $(e \in E)$ 
becomes the degenerate curve over $A_{0}$ with dual graph $\Delta$ which is 
obtained from the collection of $P_{v} := {\mathbb P}^{1}_{A_{0}}$ $(v \in V)$ 
by identifying the points 
$\alpha_{e} \in P_{v_{e}}$ and $\alpha_{-e} \in P_{v_{-e}}$ ($e \in E$), 
where $\alpha_{h}$ denotes $\infty$ if $h \in {\cal E}_{\infty}$. 

\item 
$C_{\Delta}$ gives a universal deformation 
of degenerate curves with dual graph $\Delta$, i.e. 
if $R$ is a noetherian and normal complete local ring with residue field $k$, 
and $C$ is a stable curve over $R$ with nonsingular generic fiber 
such that the closed fiber $C \otimes_{R} k$ is a degenerate curve 
with dual graph $\Delta$, in which all double points are $k$-rational, 
then there exists a ring homomorphism $A \rightarrow R$ 
giving $C_{\Delta} \otimes_{A} R \cong C$.  

\item 
$C_{\Delta} \otimes_{A} B$ is smooth over $B$ 
and is Mumford uniformized (cf. \cite{Mu}) by $\Gamma$. 

\item 
Let $\alpha_{h}$ $(h \in {\cal E})$ be complex numbers 
such that $\alpha_{e} \neq \alpha_{-e}$ 
and that $\alpha_{h} \neq \alpha_{h'}$ if $h \neq h'$ and $v_{h} = v_{h'}$. 
Then for nonzero complex numbers $q_{e}$ $(e \in E)$ 
with sufficiently small absolute value, 
$C_{\Delta}$ becomes a Riemann surface which is Schottky uniformized (cf. \cite{S}) 
by the Schottky group $\Gamma$ over ${\mathbb C}$. 

\end{itemize}

We apply the above result to the construction of a uniformized deformation of 
a degenerate pointed curve 
(see \cite[Section 2, Theorems 1 and 10]{IN} when the degenerate pointed curve 
consists of smooth pointed projective lines). 
Let $\Delta = (V, E, T)$ be a stable graph with numbering $\nu$ of $T$. 
We define its extension $\tilde{\Delta} = ( \tilde{V}, \tilde{E} )$ 
as a stable graph without tail by adding a vertex with a loop to the end 
distinct from $v_{h}$ for each tail $h \in T$. 
Then from the uniformized curve associated with $\tilde{\Delta}$, 
by substituting $0$ for the deformation parameters which correspond to 
$e \in \tilde{E} - E$ and by replacing the singular projective lines 
which correspond to $v \in \tilde{V} - V$ with marked points, 
one has the required universal deformation. 
\vspace{4ex}

\noindent
{\bf 3. Teichm\"{u}ller groupoid} 
\vspace{2ex}

\noindent
{\it 3.1. Moduli space of pointed curves} 

\noindent
We review fundamental facts on the moduli space of pointed curves 
and its compactification \cite{DM, KM, K}. 
Let $g$ and $n$ be nonnegative integers such that $n$ and $2g - 2 + n$ are positive. 
As in the introduction, 
${\cal M}_{g, n}$ (resp. ${\cal M}_{g, \vec{n}}$) denote the moduli stacks 
over ${\mathbb Z}$ classifying proper smooth curves of genus $g$ 
with $n$ marked points (resp. $n$ marked points with nonzero tangent vectors). 
Then ${\cal M}_{g, \vec{n}}$ becomes naturally a principal $({\mathbb G}_{m})^{n}$-bundle 
on ${\cal M}_{g, n}$. 
Furthermore, let $\overline{\cal M}_{g, n}$ denote the Deligne-Mumford-Knudsen 
compactification of ${\cal M}_{g, n}$ which is defined as the moduli stack 
over ${\mathbb Z}$ classifying stable curves of genus $g$ with $n$ marked points, 
and $\overline{\cal M}_{g, \vec{n}}$ denote the $({\mathbb A}^{1})^{n}$-bundle 
on $\overline{\cal M}_{g, n}$ obtained by considering (possibly zero) tangent vectors 
on marked points. 
For these moduli stacks ${\cal M}_{*,*}$ and $\overline{\cal M}_{*,*}$, 
${\cal M}_{*,*}^{\rm an}$ and $\overline{\cal M}_{*,*}^{\rm an}$ denote 
the associated complex orbifolds. 
A {\it point at infinity} on ${\cal M}_{g, n}$ (resp. ${\cal M}_{g, \vec{n}}$) is a point on 
$\overline{\cal M}_{g, n}$ (resp. $\overline{\cal M}_{g, \vec{n}}$) which corresponds to 
a maximally degenerate $n$-pointed curve, 
and a {\it tangential point at infinity} is a point at infinity with tangential structure 
over ${\mathbb Z}$. 
\vspace{2ex}

\noindent
{\it 3.2. Coordinates on moduli space} 

\noindent
To obtain explicit local coordinates on the above moduli stacks, 
we will rigidify a coordinate on each projective line appearing 
as an irreducible component of the base degenerate curve. 
In the maximally degenerate case, this process is considered in \cite{IN} 
using the notion of ``tangential structure''. 
A rigidification of an oriented stable graph $\Delta = (V, E, T)$ with numbering
$\nu$ of $T$ means a collection $\tau = \left( \tau_{v} \right)_{v \in V}$ 
of injective maps
$$
\tau_{v} : \{ 0, 1, \infty \} \rightarrow 
\left\{ h \in \pm E \cup T \ | \ v_{h} = v \right\} 
$$
such that $\tau_{v}(a) \neq - \tau_{v'}(a)$ for any $a \in \{ 0, 1, \infty \}$ 
and distinct elements $v, v' \in V$ with $\tau_{v}(a), \tau_{v'}(a) \in \pm E$. 
One can see that any stable graph has a rigidification by the induction
on the number of edges and tails. 
Denote by $A_{(\Delta, \tau)}$ the ring of formal power series ring of $q_{e}$ $(e \in E)$  over the ${\mathbb Z}$-algebra which is generated by 
$\alpha_{h}$ ($h \in {\cal E}$), $1/(\alpha_{e} - \alpha_{-e})$ ($e, -e \in {\cal E} - T$) 
and $1/(\alpha_{h} - \alpha_{h'})$ ($h, h' \in {\cal E}$ with $h \neq h'$ and $v_{h} = v_{h'}$), 
where  
$$ 
{\cal E} = \pm E \cup T - \left\{ \tau_{v}(\infty) \ | \ v \in V \right\}, 
$$
$\alpha_{h} = a$ for $h = \tau_{v}(a)$ $(a \in \{ 0, 1 \})$ and $\alpha_{h}$ are variables. 
Then as is stated above (see also \cite[2.3.9]{IN}), 
there exists a stable $\sharp T$-pointed curve $C_{(\Delta, \tau)}$ over 
$A_{(\Delta, \tau)}$ which is obtained as the quotient by $\pi_{1}(\Delta)$ 
of the glued scheme of pointed projective lines associated with 
the universal cover of $\Delta$. 
Therefore, $C_{(\Delta, \tau)}$ gives a universal deformation by $q_{e}$ $(e \in E)$ of 
the degenerate $\sharp T$-pointed curve with dual graph $\Delta$. 

Note that if one takes another rigidification, 
then $q_{e}$ $(e \in E)$ may give different deformation parameters of 
the degenerate pointed curve, 
and these parameters associated with distinct rigidifications can be compared 
(cf. \cite[Section 2]{I3}). 

Let $\tau$ be a rigidification of an oriented stable graph 
$\Delta = (V, E, T)$ with numbering of $T$, 
and put
$$
{\cal E}_{\tau} = \pm E \cup T - \bigcup_{v \in V} {\rm Im}(\tau_{v}). 
$$
Then $\alpha_{h}$ $(h \in {\cal E}_{\tau})$ and $q_{e}$ $(e \in E)$ give effective parameters 
of the moduli and the deformation of degenerate $\sharp T$-pointed curves 
with dual graph $\Delta$ respectively. 
Therefore, if we put $g = {\rm rank}_{\mathbb Z} H_{1} (\Delta, {\mathbb Z})$ 
and $n = \sharp T$, 
then
$$
\left( \alpha_{h} \ (h \in {\cal E}_{\tau}), q_{e} \ (e \in E) \right) 
$$
gives a system of formal coordinates on an \'{e}tale neighborhood of $Z_{\Delta}$, 
where $Z_{\Delta}$ denotes the substack of $\overline{\cal M}_{g, n}$ 
classifying  degenerate $n$-pointed curves with dual graph $\Delta$. 
Furthermore, by the above result, 
this system gives local coordinates on an \'{e}tale neighborhood of 
the complex orbifold $Z_{\Delta}^{\rm an}$ associated with $Z_{\Delta}$. 
\vspace{2ex}

\noindent
{\it 3.3. Teichm\"{u}ller groupoid} 

\noindent
Let ${\cal M}_{g, \vec{n}}^{\rm an}$ denote the complex orbifold associated with 
${\cal M}_{g, \vec{n}}$ as above. 
Then the {\it Teichm\"{u}ller groupoid} for ${\cal M}_{g, \vec{n}}$ is defined as 
the fundamental groupoid for ${\cal M}_{g, \vec{n}}^{\rm an}$ with tangential base points 
which correspond to maximally degenerate pointed curves. 
Its fundamental paths called {\it basic moves} are half-Dehn twists, 
fusing moves and simple moves defined as follows. 

Let $\Delta = (V, E, T)$ be a trivalent graph as above, 
and assume that $\Delta$ is trivalent. 
Then for any rigidification $\tau$ of $\Delta$, 
$\pm E \cup T = \bigcup_{v \in V} {\rm Im}(\tau_{v})$, 
and hence $A_{(\Delta, \tau)}$ is the formal power series ring over ${\mathbb Z}$ 
of $3g - 3 + n$ variables $q_{e}$ $(e \in E)$. 
First, the {\it half-Dehn twist} $\delta_{e}^{1/2}$ associated with $e$ is defined as 
the deformation of the pointed Riemann surface corresponding to $C_{\Delta}$ 
by $q_{e} \mapsto - q_{e}$. 
Second, a {\it fusing move (or associative move, A-move)} is defined to be 
different degeneration processes of a $4$-hold Riemann sphere. 
A fusing move changes $(\Delta, e)$ to another trivalent graph $(\Delta', e')$ 
such that $\Delta$, $\Delta'$ become the same graph, 
which we denote by $\Delta''$, if $e, e'$ shrink to a point. 
We denote this move by $\varphi(e, e')$. 
In \cite[Section 2]{I3}, this move is constructed using $C_{\Delta''}$. 
Finally,  {\it simple move (or S-move)} is defined to be different degeneration processes 
of a $1$-hold complex torus. 

Then as the {\it completeness theorem} called in \cite{MS}, 
the following assertion is conjectured in \cite{Gr2} and 
shown in \cite{BK1, BK2, FG, HLS, MS, NS} (especially in \cite[Sections 7 and 8]{NS}  using the notion of quilt-decompositions of Riemann surfaces). 
\vspace{2ex}

\noindent
COMPLETENESS THEOREM 
\begin{it}

\noindent
The Teichm\"{u}ller groupoid is generated by half-Dehn twists, 
fusing moves and  simple moves with relations induced from 
${\cal M}_{0, \vec{4}}$, ${\cal M}_{0, \vec{5}}$, ${\cal M}_{1, \vec{1}}$ 
and ${\cal M}_{1, \vec{2}}$. 
\end{it}
\vspace{4ex}

\noindent
{\bf 4. Mixed Teichm\"{u}ller motives} 
\vspace{2ex}

\noindent
{\it 4.1. Definition of mixed Teichm\"{u}ller motives} 

\noindent
Let ${\mathbb H} = \left( {\mathbb H}^{\rm Be}, {\mathbb H}^{\rm dR} \right)$ 
denote the (pure) motivic local system on ${\cal M}_{g, \vec{n}}$ 
associated with $R^{1} \pi_{*} {\mathbb Q}$ for the universal proper smooth curve 
$\pi : {\cal C} \rightarrow {\cal M}_{g, \vec{n}}$. 
Namely, ${\mathbb H}^{\rm Be}$ is the ${\mathbb Q}$-local system 
on ${\cal M}_{g, \vec{n}}^{\rm an}$ given by 
${\mathbb H}^{\rm Be}_{s} = H_{\rm Be}^{1} \left( {\cal C}_{s}, {\mathbb Q} \right)$   
$\left( s \in {\cal M}_{g, \vec{n}}^{\rm an} \right)$, 
and ${\mathbb H}^{\rm dR}$ is the vector bundle 
on ${\cal M}_{g, \vec{n}/{\mathbb Q}} = {\cal M}_{g, \vec{n}} \otimes {\mathbb Q}$ 
with Gauss-Manin connection given by 
${\mathbb H}^{\rm dR}_{s} = H_{\rm dR}^{1} \left( {\cal C}_{s} \right)$   
$\left( s \in {\cal M}_{g, \vec{n}/{\mathbb Q}} \right)$ with canonical isomorphism 
${\mathbb H}^{\rm Be}_{s} \otimes {\mathbb C} \cong 
{\mathbb H}^{\rm dR}_{s} \otimes {\mathbb C}$.
\vspace{2ex}

\noindent
{\it Definition.} 

\noindent 
A {\it mixed Teichm\"{u}ller motive} ${\mathbb V}$ (over ${\mathbb Z}$) 
for ${\cal M}_{g, \vec{n}}$ is defined as a unipotent motivic local system 
$$
{\mathbb V} = \left( {\mathbb V}^{\rm Be}, \ {\mathbb V}^{\rm dR} \right),
$$
namely 
\begin{itemize}

\item[{\rm (i)}] 
${\mathbb V}^{\rm Be}$ is a ${\mathbb Q}$-local system on 
${\cal M}_{g, \vec{n}}^{\rm an}$, 

\item[{\rm (ii)}]  
${\mathbb V}^{\rm dR}$ is a filtered vector bundle ${\cal V}$ 
on  $\overline{\cal M}_{g, \vec{n}/{\mathbb Q}}$ 
with Hodge filtration $F^{\bullet}$ and flat connection 
$$
\nabla : {\cal V} \rightarrow {\cal V} \otimes 
\Omega^{1}_{\overline{\cal M}_{g, \vec{n}/{\mathbb Q}}} 
\left( \log \left( {\cal D}_{g, \vec{n}/{\mathbb Q}} \right) \right); \ 
{\cal D}_{g, \vec{n}} = \overline{\cal M}_{g, \vec{n}} - {\cal M}_{g, \vec{n}}
$$
which has nilpotent residue along each component of 
${\cal D}_{g, \vec{n}/{\mathbb Q}}$ and satisfies the Griffith transversality: 
$$
\nabla : F^{p} {\cal V} \rightarrow F^{p-1} {\cal V} \otimes 
\Omega^{1}_{\overline{\cal M}_{g, \vec{n}/{\mathbb Q}}} 
\left( \log \left( {\cal D}_{g, \vec{n}/{\mathbb Q}} \right) \right) 
$$ 

\end{itemize}
such that  
\begin{itemize}

\item[{\rm (iii)}] 
there exists an increasing weight filtration $W_{\bullet}$ of ${\mathbb V}$ such that 
$$
\bigcup_{m} W_{m} {\mathbb V} = {\mathbb V}, \ \ 
\bigcap_{m} W_{m} {\mathbb V} = \{ 0 \} 
$$
and that each weight graded quotient ${\rm Gr}^{W}_{m} {\mathbb V}$ 
of ${\mathbb V}$ is isomorphic to a direct sum of subquotients of copies 
(with multiplicities) of ${\mathbb H}^{\otimes (m+2r)}(r)$, 

\item[{\rm (iv)}]  
there exists a comparison isomorphism between 
$\left( {\mathbb V}^{\rm dR}, W_{\bullet} \right) \otimes 
{\cal O}_{\overline{\cal M}_{g, \vec{n}}^{\rm an}}$ 
with the canonical extension of the filtered flat bundle 
$\left( {\mathbb V}^{\rm Be}, W_{\bullet} \right) \otimes 
{\cal O}_{{\cal M}_{g, \vec{n}}^{\rm an}}$ 
to $\overline{\cal M}_{g, \vec{n}}^{\rm an}$.

\end{itemize}
Furthermore, ${\mathbb V}$ is required to have the structure of 
mixed Tate motives over ${\mathbb Z}$ on points at infinity as follows: 
\begin{itemize}

\item[{\rm (v)}] 
for each tangential point $t$ at infinity on ${\cal M}_{g, \vec{n}}$, 
there exists a mixed Tate motive $V(t)$ over ${\mathbb Z}$ satisfying 
\begin{itemize}

\item[$\bullet$] 
the fiber ${\mathbb V}^{\rm dR}_{t}$ over $t$ gives the de Rham realization $V(t)^{\rm dR}$ 
of $V(t)$, 

\item[$\bullet$] 
the monodromy representation associated with ${\mathbb V}^{\rm Be}$ 
gives a homomorphism
$\pi_{1} \left( {\cal M}_{g, \vec{n}}^{\rm an}; t \right) \rightarrow 
{\rm Aut} \left( V(t)^{\rm Be} \right)$, 
where $V(t)^{\rm Be}$ denotes the Betti realization of $V(t)$, 

\item[$\bullet$] 
the fiber of the isomorphism in (iv) over $t$ gives the comparison isomorphism 
$\left( V(t)^{\rm dR}, W_{\bullet} \right) \otimes {\mathbb C} \cong 
\left( V(t)^{\rm Be}, W_{\bullet} \right) \otimes {\mathbb C}$ 
for the motive $V(t)$,  

\end{itemize} 

\item[{\rm (vi)}] 
for each prime number $l$ and tangential points $t_{1}, t_{2}$ at infinity on 
${\cal M}_{g, \vec{n}}$, 
the monodromy isomorphism $V(t_{1})^{\rm Be} \stackrel{\sim}{\rightarrow} V(t_{2})^{\rm Be}$  along a path in ${\cal M}_{g, \vec{n}}^{\rm an}$ from $t_{1}$ to $t_{2}$ gives rise to 
a ${\rm Gal} \left( \overline{\mathbb Q}/{\mathbb Q} \right)$-equivariant isomorphism 
$V(t_{1})^{l} \stackrel{\sim}{\rightarrow} V(t_{2})^{l}$, 
where $V(t)^{l}$ denotes the $l$-adic realization of $V(t)$. 

\end{itemize}

The category of mixed Teichm\"{u}ller motives for ${\cal M}_{g, \vec{n}}$ is 
a ${\mathbb Q}$-linear neutral tannakian category 
which we denote by ${\sf MTeM}_{g, \vec{n}}$.  

A {\it mixed Tate local system} on $({\mathbb G}_{m})^{d}$ is defined as 
a unipotent motivic local system on $({\mathbb G}_{m})^{d}$ which satisfy 
the above conditions (i)--(v) 
replacing ${\cal M}_{g, \vec{n}}$, $\overline{\cal M}_{g, \vec{n}}$ 
by $({\mathbb G}_{m})^{d}$, $({\mathbb A}^{1})^{d}$ respectively. 
A mixed Teichm\"{u}ller motive for ${\cal M}_{g, \vec{n}}$ 
gives a mixed Tate local system on $({\mathbb G}_{m})^{3g-3+2n}$ 
around each point at infinity. 
\vspace{2ex}
  
\noindent
PROPOSITION 4.1 

\begin{it} 
\noindent
The category of mixed Tate local systems on $({\mathbb G}_{m})^{d}$ is 
equivalent to the category of mixed Tate motives over ${\mathbb Z}$. 
\end{it}
\vspace{2ex} 

\noindent
{\it Proof} 

\noindent
First, assume that $d = 1$, and let ${\mathbb V}$ be a mixed Tate local system on 
${\mathbb G}_{m} = {\mathbb P}^{1} - \{ 0, \infty \}$. 
Extend ${\mathbb V}$ to a unipotent ${\mathbb Q}$-local system on 
${\mathbb P}^{1} - \{ 0, 1, \infty \}$ with trivial monodromy around $1$ 
which we denote by the same symbol. 
Then by a result of Brown \cite{Br}, 
one may assume that the fiber ${\mathbb V}_{t}$ of ${\mathbb V}$ 
over $t = \overrightarrow{01}$ is derived from the motivic fundamental group of  ${\mathbb P}^{1} - \{ 0, 1, \infty \}$. 
Since ${\mathbb V}$ gives a torsor of ${\mathbb V}_{t}$ by a closed path around $0$, 
it is an extension of ${\mathbb Q}(0)$ by ${\mathbb V}_{t}$ 
as a unipotent ${\mathbb Q}$-local system on ${\mathbb P}^{1} - \{ 0, 1, \infty \}$. 
Therefore, ${\mathbb V}$ becomes the Betti trivialization of 
a mixed Tate motive over ${\mathbb Z}$, 
and hence the assertion for any $d$ follows from induction. 
\ $\square$ 
\vspace{2ex}

\noindent
{\it 4.2. Monodromy for mixed Teichm\"{u}ller motives} 

\noindent
Let $\Delta$ be a trivalent graph of type $(g, n)$, 
namely satisfying ${\rm rank}_{\mathbb Z} H_{1}(\Delta, {\mathbb Z}) = g$ and having 
a bijection $\{ \mbox{tails of $\Delta$} \} \stackrel{\sim}{\rightarrow} \{ 1,..., n \}$. 
Let $\Delta'$ be a trivalent subtree of $\Delta$ with $m$ tails 
with numbering from $1$ to $m$, 
and denote by $C'$ the (possibly nonconnected) maximally degenerate pointed curve 
with dual graph $\Delta - (\Delta' \ \mbox{without tails})$. 
Then by attaching $C'$ to stable curves over $\overline{\cal M}_{0, m}$ 
under this numbering, 
one has a morphism 
$$
\iota_{\Delta, \Delta'} : \overline{\cal M}_{0, m} \rightarrow 
\overline{\cal M}_{g, \vec{n}}.
$$

\noindent
PROPOSITION 4.2 

\begin{it} 
\noindent
Let ${\mathbb V}$ be a mixed Teichm\"{u}ller motive for ${\cal M}_{g, \vec{n}}$, 
and denote by ${\cal V}$ the associated vector bundle on 
$\overline{\cal M}_{g, \vec{n}/{\mathbb Q}}$. 
Then the pullback $\left( \iota_{\Delta, \Delta'} \right)^{*}({\cal V})$ of ${\cal V}$ 
by $\iota_{\Delta, \Delta'}$ becomes a trivial bundle on $\overline{\cal M}_{0,m}$, 
and this trivialization is unique up to a change of basis of a fiber of ${\cal V}$. 
Furthermore, if ${\mathbb W}$ is a mixed Teichm\"{u}ller submotive of ${\mathbb V}$, 
then for the vector bundle ${\cal W}$ associated with ${\mathbb W}$, 
the trivialization of $\left( \iota_{\Delta, \Delta'} \right)^{*}({\cal W})$ 
is compatible with that of $\left( \iota_{\Delta, \Delta'} \right)^{*}({\cal V})$. 
\end{it}
\vspace{2ex}

\noindent
{\it Proof} 

\noindent
Denote the canonical extension of ${\mathbb H}$ to $\overline{\cal M}_{g, \vec{n}}$ 
by the same symbol. 
Then $\left( \iota_{\Delta, \Delta'} \right)^{*}({\mathbb H})$ has the weight filtration: 
$$
0 \rightarrow {\mathbb Q}^{\oplus (2g+n-m)} \rightarrow 
\left( \iota_{\Delta, \Delta'} \right)^{*}({\mathbb H}) \rightarrow  
{\mathbb Q}(-1)^{\oplus (m-1)} \rightarrow 0. 
$$
Since $\overline{\cal M}_{0,m}$ is rational, 
$H^{1} \left( \overline{\cal M}_{0,m}, {\cal O}_{\overline{\cal M}_{0,m}} \right) = \{ 0 \}$, 
and hence $\left( \iota_{\Delta, \Delta'} \right)^{*}({\mathbb H})$ is trivialized. 
Therefore, by the above condition (iii), 
$\left( \iota_{\Delta, \Delta'} \right)^{*}({\cal V})$ is also trivialized. 
The uniqueness follows from that any morphism from the complete variety $\overline{\cal M}_{0,m}$ to the affine variety $GL_{r}$ $(r = {\rm rank}({\cal V}) )$ 
becomes a constant map. 

Let ${\cal W}$ be as above. 
Then $\left( \iota_{\Delta, \Delta'} \right)^{*}({\cal W})$ and 
$\left( \iota_{\Delta, \Delta'} \right)^{*}({\cal V}/{\cal W})$ are trivialized, 
and hence by the vanishing of 
$H^{1} \left( \overline{\cal M}_{0,m}, {\cal O}_{\overline{\cal M}_{0,m}} \right)$, 
$\left( \iota_{\Delta, \Delta'} \right)^{*}({\cal V})$ is uniquely trivialized 
in a compatible way with the above trivializations. 
\ $\square$
\vspace{2ex}

\noindent
COROLLARY 4.3 

\begin{it} 
\noindent
Let ${\mathbb V}$ and ${\cal V}$ be as above. 
Then ${\cal V}$ can be trivialized along all fusing moves in a unique way 
up to a change of basis of a fiber of ${\cal V}$. 
Furthermore, if ${\mathbb W}$ is a mixed Teichm\"{u}ller submotive of ${\mathbb V}$, 
then for the vector bundle ${\cal W}$ associated with ${\mathbb W}$, 
the trivialization of ${\cal W}$ is compatible with that of ${\cal V}$. 
\end{it}
\vspace{2ex}

\noindent
{\it Proof} 

\noindent
By the completeness theorem, 
the points at infinity on ${\cal M}_{g, \vec{n}}^{\rm an}$ are connected 
by compositions of fusing moves and half-Dehn twists with relations 
induced from ${\cal M}_{0, \vec{4}}^{\rm an}$ and ${\cal M}_{0, \vec{5}}^{\rm an}$. 
Then the assertion follows from Proposition 4.2 
\ $\square$
\vspace{2ex}

For elements $A, B$ of ${\rm End}_{\mathbb C} \left( {\mathbb C}^{r} \right)$, 
We consider the differential equation 
$$
G'(t) = \left( \frac{A}{t} + \frac{B}{t-1} \right) G(t) \ \ (0 < t < 1),  
$$
and take its ${\rm End}_{\mathbb C} \left( {\mathbb C}^{r} \right)$-valued solutions 
$G_{i}(t)$ $(i = 0, 1)$ which are normalized in the sense that 
$$
\lim_{t \rightarrow 0} \frac{G_{0}(t)}{t^{A}} = 
\lim_{t \rightarrow 1} \frac{G_{1}(t)}{(1-t)^{B}} = 1. 
$$ 
Then the connection matrix $\Phi (A, B)$ is defined as $G_{1}(t)^{-1} \cdot G_{0}(t)$. 
\vspace{2ex}

\noindent
THEOREM 4.4 
 
\begin{it}
\noindent 
Denote by $\delta_{e}^{1/2}$ and $\varphi(e, e')$ the half Dehn twist and 
the fusing move respectively given in 3.3. 
For a mixed Teichm\"{u}ller motive ${\mathbb V}$, 
let ${\cal V}$ be the associated vector bundle with connection $\nabla$ on 
$\overline{\cal M}_{g, \vec{n}/{\mathbb Q}}$.  
Then under a trivialization of ${\cal V}$ by ${\mathbb C}^{r}$ along $\varphi(e, e')$,  
the monodromy of $\nabla$ along $\delta_{e}^{1/2}$ and $\varphi(e, e')$ are expressed as $\exp \left( \pi \sqrt{-1} \cdot {\rm Res}_{e} (\nabla) \right)$ and  
$\Phi \left( {\rm Res}_{e} (\nabla), {\rm Res}_{e'} (\nabla) \right)$ respectively, 
where ${\rm Res}_{e} (\nabla)$ denotes the residue of $\nabla$ around $q_{e} = 0$. 
\end{it}
\vspace{2ex}

\noindent
{\it Proof} 

\noindent
This assertion follows from \cite[Proposition 1 and Theorem 2]{I4}. 
\ $\square$
\vspace{2ex}

\noindent
{\it 4.3. Category of mixed Teichm\"{u}ller motives} 

\noindent
Let ${\cal T}_{g, n}$ be the set of trivalent graphs of type $(g, n)$ 
which is identified with the set of points at infinity on ${\cal M}_{g, n}$, 
and ${\cal L}_{g, n}$ be the set of loops in all trivalent graphs of type $(g, n)$. 
Then there is a natural map ${\cal L}_{g, n} \rightarrow {\cal T}_{g, n}$ 
which sends each $\ell \in {\cal L}_{g, n}$ to the element $t_{\ell}$ of ${\cal T}_{g, n}$ 
containing $\ell$. 
Let ${\mathbb V}$ be a mixed Teichm\"{u}ller motive for ${\cal M}_{g, \vec{n}}$. 
Then for each $t \in {\cal T}_{g, n}$, 
${\mathbb V}_{t}$ denotes the mixed Tate motive over ${\mathbb Z}$ 
given by the fiber of ${\mathbb V}$ at the point at infinity corresponding to $t$. 
By the definition of universal mixed elliptic motives \cite{HM}, 
for each $\ell \in {\cal L}_{g, n}$, 
the restriction of ${\mathbb V}$ to the space ${\cal M}_{1, \vec{1}}$ corresponding to 
$\ell$ gives rise to a universal mixed elliptic motive 
which we denote by ${\mathbb V}_{\ell}$. 

We define the {\it $(g, n)$-type product category} ${\sf PTE}_{g, \vec{n}}$ 
of mixed Tate motives over ${\mathbb Z}$ and universal mixed elliptic motives 
as follows:   
\begin{itemize}

\item 
Each object consists of mixed Tate local systems ${\mathbb V}(t)$ on 
$({\mathbb G}_{m})^{3g-3+2n}$ indexed by $t \in {\cal T}_{g,n}$ 
and universal mixed elliptic motives ${\mathbb V}(\ell)$ indexed by $\ell \in {\cal L}_{g, n}$ 
with identifications ${\mathbb V}(t)^{\rm dR} \cong V$ as ${\mathbb Q}$-vector spaces 
for a fixed ${\mathbb Q}$-vector space $V$. 
Furthermore, for each $\ell \in {\cal L}_{g, n}$, 
the fiber of ${\mathbb V}(\ell)$ over the point at infinity on ${\cal M}_{1, \vec{1}}$ 
is identified with the restriction of ${\mathbb V}(t_{\ell})$. 

\item 
Each morphism 
$\left( {\mathbb V}(t), {\mathbb V}(l) \right) \rightarrow 
\left( {\mathbb V}'(t), {\mathbb V}'(l) \right)$ 
consists of morphisms 
${\mathbb V}(t) \rightarrow {\mathbb V}'(t)$ $\left( t \in {\cal T}_{g, n} \right)$ 
as mixed Tate motives over ${\mathbb Z}$ and those 
${\mathbb V}(l) \rightarrow {\mathbb V}'(l)$ $\left( \ell \in {\cal L}_{g, n} \right)$ 
as universal mixed elliptic motives 
which are compatible with the above structure of these objects. 

\end{itemize} 

\noindent
THEOREM 4.5
 
\begin{it} 
\noindent
By the functor 
$$
{\mathbb V} \mapsto 
\left( {\mathbb V}_{t} \ (t \in {\cal T}_{g, n}), 
{\mathbb V}_{\ell} \ (\ell \in {\cal L}_{g, n}) \right), 
$$
${\sf MTeM}_{g, \vec{n}}$ is equivalent to a full subcategory of ${\sf PTE}_{g, \vec{n}}$. 
\end{it}
\vspace{2ex}

\noindent
{\it Proof} 

\noindent
By Corollary 4.3, 
$\left( {\mathbb V}_{t}, {\mathbb V}_{\ell} \right)$ 
becomes an object of ${\sf PTE}_{g, \vec{n}}$. 
Therefore, to prove the assertion, 
it is enough to show that for any object 
${\mathbb V} = 
\left( {\mathbb V}^{\rm Be}, {\mathbb V}^{\rm dR} = ({\cal V}, \nabla) \right)$ 
of ${\sf MTeM}_{g, \vec{n}}$, 
$$
{\rm Hom}_{{\sf MTeM}_{g, \vec{n}}} ({\mathbb Q}(0), {\mathbb V}) \rightarrow 
{\rm Hom}_{{\sf PTE}_{g, \vec{n}}} ({\mathbb Q}(0), {\mathbb V}) 
$$
is an isomorphism. 
Let $\psi : {\mathbb Q}(0) \rightarrow {\mathbb V}$ be a nonzero homomorphism 
in ${\sf PTE}_{g, \vec{n}}$. 
Then for the trivialization $V$ of ${\cal V}$, 
the subspace of $V$ corresponding to $\psi$ is stable under the action of 
${\rm Res}_{e}(\nabla)$, 
and hence by Theorem 4.4, 
$\psi$ gives a $1$-dimensional trivial subbundle of ${\mathbb V}^{\rm dR}$ 
along all fusing moves. 
By definition, 
$\psi$ also corresponds to a $1$-dimensional trivial subbundle of 
${\mathbb V}^{\rm dR}$ along each simple move. 
Therefore, by the completeness theorem, 
there exists uniquely a $1$-dimensional trivial subbundle of ${\mathbb V}^{\rm dR}$ 
corresponding to $\psi$ on $\overline{\cal M}_{g, \vec{n}}$ 
which gives nonzero homomorphism ${\mathbb Q}(0) \rightarrow {\mathbb V}$ 
in ${\sf MTeM}_{g, \vec{n}}$. 
\ $\square$
\vspace{2ex}

The ${\mathbb Q}$-linear neutral tannakian categories 
${\sf MTeM}_{g, \vec{n}}$ and ${\sf PTE}_{g, \vec{n}}$ are identified with 
the categories of representations of affine group schemes 
which we will denote by $\pi_{1} \left( {\sf MTeM}_{g, \vec{n}} \right)$ and 
$\pi_{1} \left( {\sf PTE}_{g, \vec{n}} \right)$ respectively, 
and call the motivic  Galois groups. 
\vspace{2ex}

\noindent
COROLLARY 4.6 

\begin{it} 
\noindent
There exists a natural surjective homomorphism 
\end{it}
$$
\pi_{1} \left( {\sf PTE}_{g, \vec{n}} \right) \rightarrow 
\pi_{1} \left( {\sf MTeM}_{g, \vec{n}} \right). 
$$

\noindent
{\it Proof} 

\noindent
The assertion follows from Theorem 4.5 which states that 
${\sf MTeM}_{g, \vec{n}}$ is a full subcategory of ${\sf PTE}_{g, \vec{n}}$. 
\ $\square$
\vspace{4ex}

\noindent
{\bf 5. Motivic fundamental torsors}
\vspace{2ex}

\noindent
{\it 5.1. Fundamental torsors of curves} 

\noindent
Assume that $n$ is a positive integer, 
and let $C^{\circ}$ be an algebraic curve over a subfield $K$ of ${\mathbb C}$ which is  obtained from a proper smooth curve $C$ of genus $g$ by removing $n$ points. 
Note that $C^{\circ}$ is not complete, and hence its first de Rham cohomology group 
$H^{1}_{\rm dR}(C^{\circ})$ has a basis $B_{C^{\circ}}$ consisting of 
$2g + n-1$ meromorphic $1$-forms on $C$ of the first or second kind 
which may have poles outside $C^{\circ}$. 
For $w_{1},..., w_{m} \in B_{C^{\circ}}$, 
we define a ${\cal D}$-module $D(w_{1},..., w_{m})$ on $C^{\circ}$ 
whose underlying bundle is given by the trivial bundle $K^{m+1} \times C^{\circ}$ 
with connection form $d - \sum_{i = 1}^{m} e_{i, i+1} w_{i}$, 
where $e_{i,j}$ denotes the square matrix of degree $m+1$ whose $(k, l)$-entry is 
$\delta_{ik} \cdot \delta_{jl}$ $(\delta_{ij}$ denotes Kronecker's delta). 
We consider the tannakian subcategory of ${\cal D}$-modules on $C^{\circ}$ 
generated by $D(w_{1},..., w_{m})$ $\left( w_{1},..., w_{m} \in B_{C^{\circ}} \right)$. 
Since the underlying bundles of objects in this category are trivial, 
for each $K$-rational (tangential) point $x$ on $C^{\circ}$, 
one can define the fiber functor on this category by taking the (trivial) fibers at $x$. 
Denote by $\pi_{1}^{\rm dR}(C^{\circ}; x)$ the tannakian fundamental group of 
this category which is a pro-finite algebraic group over $K$, 
and by ${\cal A}^{\rm dR}(C^{\circ}; x)$ the enveloping algebra of the Lie algebra 
${\rm Lie} \left( \pi_{1}^{\rm dR}(C^{\circ}; x) \right)$. 

Let $(C^{\circ})^{\rm an}$ be a Riemann surface associated with 
$C^{\circ} \otimes_{K} {\mathbb C}$, 
and for each (tangential) points $x, y$ on $C^{\circ}$, 
denote by $\pi_{1}((C^{\circ})^{\rm an}; x, y)$ the set of homotopy classes of paths 
from $x$ to $y$ in $(C^{\circ})^{\rm an}$. 
When $x = y$, $\pi_{1}((C^{\circ})^{\rm an}; x, y)$ becomes the fundamental group 
$\pi_{1}((C^{\circ})^{\rm an}; x)$ of $(C^{\circ})^{\rm an}$ with base point $x$. 
We consider the tannakian category of unipotent local systems on $(C^{\circ})^{\rm an}$ 
with fiber functor obtained from taking the fiber at $x$.  
Then it is shown in \cite{D} that its tannakian fundamental group 
$\pi_{1}^{\rm Be}((C^{\circ})^{\rm an}; x)$ is a pro-algebraic group over ${\mathbb Q}$, 
and the associated enveloping algebra ${\cal A}^{\rm Be}((C^{\circ})^{\rm an}; x)$ 
of ${\rm Lie} \left( \pi_{1}^{\rm Be}((C^{\circ})^{\rm an}; x) \right)$ is isomorphic to 
$$
\lim_{m \rightarrow \infty} 
\left. {\mathbb Q} \left[ \pi_{1} ((C^{\circ})^{\rm an}; x) \right] \right/ I^{m},  
$$
where $I$ denotes the augmentation ideal of the group ring 
${\mathbb Q} \left[ \pi_{1}((C^{\circ})^{\rm an}; x) \right]$. 
Since $\pi_{1} ((C^{\circ})^{\rm an}; x)$ is a free group of rank $2g + n-1$, 
${\cal A}^{\rm Be} ((C^{\circ})^{\rm an}; x)$ becomes the ring of 
noncommutative formal power series over ${\mathbb Q}$ in $2g + n-1$ variables. 
\vspace{2ex}

\noindent
PROPOSITION 5.1 

\begin{it}
\noindent
For a $K$-rational (tangential) point $x$ on $C^{\circ}$, 
there exists a canonical isomorphism 
$$
{\cal A}^{\rm dR} (C^{\circ}; x) \otimes_{K} {\mathbb C} \cong  
{\cal A}^{\rm Be} ((C^{\circ})^{\rm an}; x) \otimes_{\mathbb Q} {\mathbb C}.  
$$ 
Consequently, 
${\rm Lie} \left( \pi_{1}^{\rm dR} (C^{\circ}; x) \right) \otimes_{K} {\mathbb C}$ 
is isomorphic to 
${\rm Lie} \left( \pi_{1}^{\rm Be} ((C^{\circ})^{\rm an}; x) \right) \otimes_{\mathbb Q} 
{\mathbb C}$, 
and $\pi_{1}^{\rm dR} (C^{\circ}; x) \otimes_{K} {\mathbb C}$ is isomorphic to 
$\pi_{1}^{\rm Be} ((C^{\circ})^{\rm an}; x) \otimes_{\mathbb Q} {\mathbb C}$. 
\end{it}
\vspace{2ex}

\noindent
{\it Proof} 

\noindent
By associating local systems on $(C^{\circ})^{\rm an}$ with 
$D(w_{1},..., w_{m})$ $(w_{i} \in B_{C^{\circ}})$, 
one has a group homomorphism 
$$
\pi_{1}^{\rm Be} ((C^{\circ})^{\rm an}; x) \otimes_{\mathbb Q} {\mathbb C} 
\rightarrow 
\pi_{1}^{\rm dR} (C^{\circ}; x) \otimes_{K} {\mathbb C}  
$$
which gives a ring homomorphism 
$$
{\cal A}^{\rm Be} ((C^{\circ})^{\rm an}; x) \otimes_{\mathbb Q} {\mathbb C} 
\rightarrow 
{\cal A}^{\rm dR} (C^{\circ}; x) \otimes_{K} {\mathbb C}.
$$
Let $B_{m}((C^{\circ})^{\rm an})$ be the ${\mathbb C}$-vector space of iterated integrals spanned by 
$$
\int w_{1} \cdots w_{r} \ 
\left( w_{i} \in B_{C^{\circ}}, \ r \leq m \right),  
$$
and $H^{0} \left( B_{m}((C^{\circ})^{\rm an}); x \right)$ be the space consisting of 
elements of $B_{m}((C^{\circ})^{\rm an})$ whose restriction to 
$\left\{ \mbox{loops in $(C^{\circ})^{\rm an}$ based at $x$} \right\}$ is homotopy functional. 
Then it is shown in \cite[(5.3)]{H} that 
$H^{0} \left( B_{m}((C^{\circ})^{\rm an}); x \right)$ is the dual space of 
$\left. {\mathbb C} \left[ \pi_{1} ((C^{\circ})^{\rm an}; x) \right] \right/ I^{m+1}$, 
and that by taking leading terms of iterated integrals, 
one has an exact sequence 
$$
0 \rightarrow H^{0} \left( B_{m-1}((C^{\circ})^{\rm an}); x \right) \rightarrow 
H^{0} \left( B_{m}((C^{\circ})^{\rm an}); x \right) \rightarrow 
H^{1}_{\rm dR}((C^{\circ})^{\rm an})^{\otimes m}. 
$$ 
Since 
$H^{1}_{\rm dR}(C^{\circ}) \otimes_{K} {\mathbb C} \cong 
H^{1}_{\rm dR}((C^{\circ})^{\rm an}) 
\cong  H^{1}_{\rm Be}((C^{\circ})^{\rm an}) \otimes_{\mathbb Q} {\mathbb C}$, 
one has a ring homomorphism 
$$
{\cal A}^{\rm dR} (C^{\circ}; x) \otimes_{K} {\mathbb C} \rightarrow 
{\cal A}^{\rm Be} ((C^{\circ})^{\rm an}; x) \otimes_{\mathbb Q} {\mathbb C}
$$
which is an inverse map to the above map, 
and hence this ring homomorphism is an isomorphism. 
The remaining assertions follow from that the tannakian fundamental group $\pi_{1}^{\sharp}$ 
$(\sharp = {\rm dR}, {\rm Be})$ and its Lie algebra are given as the subsets of 
${\cal A}^{\sharp}$ which consist of grouplike elements and primitive elements in 
${\cal A}^{\sharp}$ respectively. 
\ $\square$ 
\vspace{2ex}

For $K$-rational (tangential) points $x, y$ on $C^{\circ}$, 
denote by $\pi_{1}^{\rm dR}(C^{\circ}; x, y)$ the tannakian fundamental left (resp. right) torsor over $\pi_{1}^{\rm dR}(C^{\circ}; x)$ 
(resp. $\pi_{1}^{\rm dR}(C^{\circ}; y)$). 
Then one has the associated torsors 
${\rm Lie} \left( \pi_{1}^{\rm dR}(C^{\circ}; x, y) \right)$ and 
${\cal A}^{\rm dR}(C^{\circ}; x, y)$ over 
${\rm Lie} \left( \pi_{1}^{\rm dR}(C^{\circ}; *) \right)$ and 
${\cal A}^{\rm dR}(C^{\circ}; *)$ $(* = x, y)$ respectively. 
Similarly one can define 
$$
\pi_{1}^{\rm Be}((C^{\circ})^{\rm an}; x, y), \ 
{\rm Lie} \left( \pi_{1}^{\rm Be}((C^{\circ})^{\rm an}; x, y) \right), \ 
{\cal A}^{\rm Be}((C^{\circ})^{\rm an}; x, y)
$$
as the torsors over $\pi_{1}^{\rm Be}((C^{\circ})^{\rm an}; *)$, 
${\rm Lie} \left( \pi_{1}^{\rm Be}((C^{\circ})^{\rm an}; *) \right)$, 
${\cal A}^{\rm Be}((C^{\circ})^{\rm an}; *)$ respectively. 
Especially, $\pi_{1}^{\rm Be}((C^{\circ})^{\rm an}; x, y)$ is the completion of 
the ${\mathbb Q}$-vector space 
${\mathbb Q} \left[ \pi_{1} ((C^{\circ})^{\rm an}; x, y) \right]$ 
generated by $\pi_{1}((C^{\circ})^{\rm an}; x, y)$ with respect to 
$I^{m} \left( {\mathbb Q} \left[ \pi_{1}((C^{\circ})^{\rm an}; x, y) \right] \right)$ 
$(m \geq 1)$. 
Then Proposition 5.1 implies: 
\vspace{2ex}

\noindent
PROPOSITION 5.2

\begin{it}
\noindent
For a $K$-rational (tangential) point $x, y$ on $C^{\circ}$, 
there exists a canonical isomorphism 
$\pi_{1}^{\rm dR} (C^{\circ}; x, y) \otimes_{K} {\mathbb C} \cong 
\pi_{1}^{\rm Be} ((C^{\circ})^{\rm an}; x, y) \otimes_{\mathbb Q} {\mathbb C}$
which gives canonical isomorphisms: 
\end{it}
\begin{eqnarray*}
{\rm Lie} \left( \pi_{1}^{\rm dR} (C^{\circ}; x, y) \right) \otimes_{K} {\mathbb C}
& \cong & 
{\rm Lie} \left( \pi_{1}^{\rm Be} ((C^{\circ})^{\rm an}; x, y) \right) \otimes_{\mathbb Q} {\mathbb C}, 
\\ 
{\cal A}^{\rm dR} (C^{\circ}; x, y) \otimes_{K} {\mathbb C} 
& \cong & 
{\cal A}^{\rm Be} ((C^{\circ})^{\rm an}; x, y) \otimes_{\mathbb Q} {\mathbb C}. 
\end{eqnarray*}

\noindent
{\it Remark} 

\noindent
For a generalized Tate curve, the isomorphism in Proposition 5.2 
is a higher genus version of elliptic polylogarithms \cite{BL, BrL, L, LR}, 
and is described by the iterated integrals of canonical meromorphic $1$-forms 
given in \cite{I1}. 
\vspace{2ex}

\noindent
{\it 5.2. Motivic sheaves of fundamental torsors}

\noindent
Let $\pi : {\cal C} \rightarrow {\cal M}_{g, \vec{n}}$ denote the universal curve  
with sections $\sigma_{1}, ..., \sigma_{n}$, 
$$
{\cal C}^{\circ} = {\cal C} - \bigcup_{k = 1}^{n} {\rm Im}(\sigma_{k})
$$ 
denote the associated open curve and 
$\left( {\cal C}_{s} = \pi^{-1}(s); (\sigma_{1})_{s},..., (\sigma_{n})_{s} \right)$
denote these fibers at $s \in {\cal M}_{g, \vec{n}}$. 
\vspace{2ex} 

\noindent
THEOREM 5.3
 
\begin{it}
\noindent
For $i, j \in \{ 1,..., n \}$, 
$(\sigma_{i})_{s}$ and $(\sigma_{j})_{s}$ denote tangential points on ${\cal C}^{\circ}_{s}$. 
Then 
$$
{\cal A} \left( {\cal C}^{\circ}; \sigma_{i}, \sigma_{j} \right) = 
\left( {\cal A}^{\rm Be} \left( {\cal C}^{\circ}_{s}; (\sigma_{i})_{s}, (\sigma_{j})_{s} \right), 
{\cal A}^{\rm dR} \left( {\cal C}^{\circ}_{s}; (\sigma_{i})_{s}, (\sigma_{j})_{s} \right) \right) 
$$ 
gives rise to a pro-object in the category of mixed Teichm\"{u}ller motives 
for ${\cal M}_{g, \vec{n}}$. 
\end{it}
\vspace{2ex}

\noindent
{\it Proof} 

\noindent
Let ${\cal A}^{\rm Be}$ denote the local system on ${\cal M}_{g, \vec{n}}^{\rm an}$ 
given by 
${\cal A}^{\rm Be} \left( {\cal C}^{\circ}_{s}; (\sigma_{i})_{s}, (\sigma_{j})_{s} \right)$ 
whose monodromy corresponds to the natural action of the mapping class group 
$\pi_{1} \left( {\cal M}_{g, \vec{n}}^{\rm an} \right)$ on 
$\pi_{1} \left( ({\cal C}^{\circ}_{s})^{\rm an} \right)$. 
Hence the action of Dehn twists on ${\cal A}^{\rm Be}$ is unipotent. 
Let ${\cal A}^{\rm dR}$ denote the vector bundle on ${\cal M}_{g, \vec{n}/{\mathbb Q}}$ 
given by 
${\cal A}^{\rm dR} \left( {\cal C}^{\circ}_{s}; (\sigma_{i})_{s}, (\sigma_{j})_{s} \right)$. 
Then by the comparison isomorphism given in Proposition 5.2, 
there exists the canonical extension of ${\cal A}^{\rm dR}$ to 
$\overline{\cal M}_{g, \vec{n}/{\mathbb Q}}$ for which 
the associated flat connection on ${\cal A}^{\rm dR}$ has nilpotent residue 
along each component of ${\cal D}_{g, \vec{n}/{\mathbb Q}}$. 
Following Hain \cite{H}, 
we recall the construction of the Hodge filtration $F^{p}$ and the weight filtration $W_{l}$ 
on ${\cal A} = \left( {\cal A}^{\rm Be}, {\cal A}^{\rm dR} \right)$.  
Put $D_{s} = {\cal C}_{s} - {\cal C}_{s}^{\circ}$, 
and denote by 
$E^{\bullet}_{s} = E^{\bullet} \left( {\cal C}_{s} \left( \log (D_{s}) \right) \right)$ 
the complex of $C^{\infty}$ forms on $({\cal C}_{s})^{\rm an}$ 
with logarithmic singularities along $D_{s}$. 
Then $H^{0} \left( B_{m} \left( E^{\bullet}_{s} \right) \right)$ is canonically isomorphic to 
the dual space of 
${\mathbb C} \left[ \pi_{1} \left( ({\cal C}^{\circ}_{s})^{\rm an} \right) \right] / I^{m}$, 
and the Hodge, weight filtrations on ${\cal A}^{\rm Be} \otimes {\mathbb C}$ 
(and hence on ${\cal A}^{\rm dR} \otimes {\mathbb C}$) 
are given by the filtrations  
\begin{eqnarray*}
F^{p} E^{\bullet}_{s} 
& = & 
\left\{ \mbox{forms with $\geq p$ holomorphic $1$-forms on ${\cal C}(\log(D_{s}))$} \right\}, 
\\ 
W_{l} E^{\bullet}_{s} 
& = & 
\left\{ \mbox{forms with $\leq l$ $1$-forms having logarithmic pole at $D_{s}$} \right\}  
\end{eqnarray*} 
of $E^{\bullet}_{s}$ respectively. 
Especially, $F^{p}$ is spanned by iterated integrals $\int w_{1} \cdots w_{r}$ for 
at least $p$ holomorphic $1$-forms on ${\cal C}(\log(D_{s}))$. 
If $v$ is a holomorphic tangent vector at $s \in {\cal M}_{g, \vec{n}}^{\rm an}$, 
then the associated derivative $\nabla_{v}$ for the above connection on 
${\cal A}^{\rm dR}$ satisfies 
$$
\nabla_{v} \left( \int w_{1} \cdots w_{r} \right) = 
\int \nabla_{v} \left( \int w_{1} \cdots w_{r-1} \right) w_{r} + 
\int \left( \int w_{1} \cdots w_{r-1} \right) \nabla_{v}(w_{r}), 
$$ 
and hence $\nabla_{v}$ satisfies the Griffith transversality. 
Since there are a canonical exact sequence 
$$
0 \rightarrow {\mathbb Q}(1)^{\oplus (n-1)} \rightarrow 
H^{\rm Be}_{1} \left( ({\cal C}^{\circ}_{s})^{\rm an}, {\mathbb Q} \right) 
\cong I/I^{2} \rightarrow 
H^{\rm Be}_{1} \left( ({\cal C}_{s})^{\rm an}, {\mathbb Q} \right) \rightarrow 0 
$$
and a surjection 
$\left( I/I^{2} \right)^{\otimes m} \rightarrow I^{m}/I^{m+1}$, 
each weight graded quotient of ${\cal A}$ is isomorphic to a direct sum of 
subquotients of copies of ${\mathbb H}^{\otimes (m+2r)}(r)$. 
Therefore, ${\cal A}$ gives a unipotent motivic local system on ${\cal M}_{g, \vec{n}}$.    

For a tangential point $t$ at infinity on ${\cal M}_{g, \vec{n}}$, 
$\pi_{1} \left( ({\cal C}^{\circ}_{t})^{\rm an} \right)$ is the amalgamated product of 
copies of $\Pi = \pi_{1} \left( {\mathbb P}^{1}_{\mathbb C} - \{ 0, 1, \infty \} \right)$ 
for the trivalent graph associated with $t$. 
Then the limit group structure of $\pi_{1} \left( ({\cal C}^{\circ}_{t})^{\rm an} \right)$ is 
described by the vertex groups $\Pi$ and the edge groups ${\mathbb G}_{m}$, 
and hence the limit motivic structure of 
${\cal A} = \left( {\cal A}^{\rm Be}, {\cal A}^{\rm dR} \right)$ becomes 
a mixed Tate motive over ${\mathbb Z}$. 
Let $l$ be a prime number. 
Then by the theory of algebraic fundamental groups \cite{Gr1}, 
${\cal A}^{l} = {\cal A}^{\rm Be} \otimes {\mathbb Q}_{l}$ gives a $l$-adic sheaf 
with Galois action. 
Furthermore, 
by the theorem of van Kampen on fundamental groups of Riemann surfaces 
and its arithmetic version given by Grothendieck-Murre \cite{GrM} and Nakamura \cite{N}, 
the $l$-adic geometric fundamental group $\pi_{1}^{l} \left({\cal C}^{\circ}_{t} \right)$ is 
the amalgameted product of copies of 
$\pi_{1}^{l} \left( {\mathbb P}^{1} - \{ 0, 1, \infty \} \right)$ with Galois action. 
Therefore, ${\cal A}$ gives mixed Tate motives over ${\mathbb Z}$ on points at infinity, 
and hence is a pro-object in the category of mixed Teichm\"{u}ller motives. 
\ $\square$ 
\vspace{2ex}

\noindent
{\it 5.3. Polylogarithmic motive} 

\noindent
We construct the {\it polylogarithmic motive} which is a pro-object 
in the category of universal mixed Teichm\"{u}ller motives. 
Let the notation be as in 5.2. 
Then for each $i \in \{ 1,..., n \}$, 
$$
L_{s} = 
{\rm Lie} \left( \pi_{1}^{\rm Be} \left( {\cal C}^{\circ}_{s}; (\sigma_{i})_{s} \right) \right) \ 
\left( s \in {\cal M}_{g, \vec{n}}^{\rm an} \right) 
$$
gives a ${\mathbb Q}$-local system $L$ on ${\cal M}_{g, \vec{n}}^{\rm an}$ 
whose central series filtration is  
$$
L_{s}^{1} = L_{s}, \ L_{s}^{k+1} = \left[ L_{s}^{k}, L_{s} \right]. 
$$ 
We define the logarithmic local system ${\cal L}{\rm og}$ as 
$L^{2} / \left[ L^{2}, L^{2} \right]$, 
and the polylogarithmic local system ${\cal P}{\rm ol}$ as 
$L / \left[ L^{2}, L^{2} \right]$. 
Then one has a natural exact sequence 
$$
0 \rightarrow {\cal L}{\rm og} \rightarrow {\cal P}{\rm ol} \rightarrow 
\pi^{*}({\cal H})  \rightarrow 0,  
$$
where ${\cal H}_{s} = H^{\rm Be}_{1} \left( {\cal C}^{\circ}_{s}, {\mathbb Q} \right)$. 
By Theorem 5.3, one can see the following:  
\vspace{2ex} 

\noindent
THEOREM 5.4
 
\begin{it}
\noindent
For each $i \in \{ 1,..., n \}$, 
${\cal P}{\rm ol}$ gives rise to a pro-object in the category of 
mixed Teichm\"{u}ller motives for ${\cal M}_{g, \vec{n}}$. 
\end{it}
\vspace{2ex} 

\noindent
{\it 5.4. Motivic monodromy and correlator} 

\noindent
Let $s$ be a point or tangential point at infinity on ${\cal M}_{g, \vec{n}/{\mathbb Q}}$, 
and take $i, j \in \{ 1,..., n \}$. 
Since the Betti realization of the fiber over $s$ gives a fiber functor of the category of 
mixed Teichm\"{u}ller motives for ${\cal M}_{g, \vec{n}}$, 
by Theorem 5.3, one has the associated group homomorphism 
$$
\pi_{1} \left( {\sf MTeM}_{g, \vec{n}} \right) \rightarrow 
{\rm Aut} \left( {\cal A}^{\rm Be}_{s; i, j} \right); \ 
{\cal A}^{\rm Be}_{s; i, j} = 
{\cal A}^{\rm Be} \left( {\cal C}^{\circ}_{s}; (\sigma_{i})_{s}, (\sigma_{j})_{s} \right). 
$$
We call this homomorphism a {\it motivic monodromy representation} of 
$\pi_{1} \left( {\sf MTeM}_{g, \vec{n}} \right)$. 

Assume that $s$ is a point at infinity. 
Then combining the homomorphism in Corollary 4.6, 
one has the representation 
$$
\pi_{1} \left( {\sf PTE}_{g, \vec{n}} \right) \rightarrow 
{\rm Aut} \left( {\cal A}^{\rm Be}_{s; i, j} \right) 
$$
whose description can be reduced to the case when $(g, n) = (0, 3), (1, 1)$ 
by the theorem of van Kampen on fundamental groups of Riemann surfaces. 

In \cite[1.11 and Section 10]{G}, 
Goncharov construct the motivic correlator map under the existence of 
the abelian category of mixed motives over a field. 
Using our results, one has the motivic correlator map as follows. 
Assume that $s \in {\cal M}_{g, \vec{n}}$ and $i = j$, 
and put $C = {\cal C}_{s}$ and $S^{*} = \{ (\sigma_{k})_{s} \ | \ k \neq i \}$. 
Denote by ${\cal CL}{\it ie}_{C, S^{*}}^{\vee}$ 
the cyclic envelope of the tensor algebra of 
$H_{\rm Be}^{1} \left( C^{\rm an}, {\mathbb Q} \right) \oplus {\mathbb Q}[S^{*}]$ 
modulo shuffle relations. 
Then by \cite[(1.28) and Proposition 8.5]{G}, 
the motivic monodromy representation 
$\pi_{1} \left( {\sf MTeM}_{g, \vec{n}} \right) \rightarrow 
{\rm Aut} \left( {\cal A}^{\rm Be}_{s; i, i} \right)$  
gives rise to the motivic correlator map 
$$
{\rm Cor}_{\sf MTeM} : {\cal CL}{\it ie}_{C, S^{*}}^{\vee}(1) \rightarrow 
{\cal L}{\it ie}_{g, \vec{n}}, 
$$
where ${\cal L}{\it ie}_{g, \vec{n}}$ is a Lie coalgebra defined as the dual to 
the Lie algebra of $\pi_{1} \left( {\sf MTeM}_{g, \vec{n}} \right)$.  
Goncharov's motivic correlator map substantially becomes ${\rm Cor}_{\sf MTeM}$ 
since its Hodge realization with period map give the Hodge correlator map 
$$
{\rm Cor}_{\rm Hod} : {\cal CL}{\it ie}_{C, S^{*}}^{\vee}(1) \rightarrow {\mathbb C} 
$$
constructed in \cite[3.2]{G}.

\renewcommand{\refname}{\centerline{\normalsize{\bf References}}}
\bibliographystyle{amsplain}

\begin{flushleft} 
Department of Mathematics, Graduate School of Science and Engineering, 
Saga University, Saga 840-8502, Japan; ichikawn@cc.saga-u.ac.jp 
\end{flushleft} 

\end{document}